\documentclass[11pt,leqno]{article} 
\usepackage{graphics}
\newtheorem{thm}{Theorem}

\newcommand{\beqa}{\begin{eqnarray}}
\newcommand{\eeqa}{\end{eqnarray}}

\newcommand{\pf}{\noindent {\bf Proof:} $\s$ }

\newcommand{\beq}{\begin{equation}}
\newcommand{\eeq}{\end{equation}}
\newcommand{\lbl}{\label}
\newcommand{\s}{\; \;}

\newcommand{\la}{\lambda}

\title{Non-singular solutions of $p$-Laplace  problems, allowing multiple changes of sign in the nonlinearity}

\author{
Philip Korman   \\ 
Department of Mathematical Sciences \\ 
University of Cincinnati \\ 
Cincinnati Ohio 45221-0025 \\
}

\date{}

\begin{document}

\maketitle
\begin{abstract}
\noindent 
For the $p$-Laplace Dirichlet problem (where $\varphi (t)=t|t|^{p-2}$, $p>1$)
\[
 \varphi(u'(x))'+ f(u(x))=0 \s\s \mbox{for $-1<x<1$}, \s u(-1)=u(1)=0
\]
assume that $f'(u)>(p-1)\frac{f(u)}{u}>0$  for $u>\gamma>0$, while $\int _u^{\gamma} f(t) \, dt < 0$ for all $u \in (0,\gamma)$. Then any positive solution, with $\max _{(-1,1)} u(x)=u(0)>\gamma$, is non-singular, no matter how many times $f(u)$ changes sign on $(0,\gamma)$. Uniqueness of solution follows.
 \end{abstract}

\begin{flushleft}
Key words:  Non-singular positive solutions, $p$-Laplace  problems. 
\end{flushleft}

\begin{flushleft}
AMS subject classification: 34B15.
\end{flushleft}

We consider positive solutions of
\begin{equation}
\label{3-1}
 \varphi(u'(x))'+ f(u(x))=0 \s\s \mbox{for $-1<x<1$}, \s u(-1)=u(1)=0 \,,
\end{equation}
where $\varphi (t)=t|t|^{p-2}$, $p>1$, so that $\varphi' (t)=(p-1)|t|^{p-2}$.
The linearized problem is
\begin{eqnarray}
\label{3-5}
& \left( \varphi'(u'(x))w'(x) \right)'+ f'(u(x))w(x)=0 \s\s \mbox{for $-1<x<1$}, \\
&  w(-1)=w(1)=0 \,. \nonumber
\end{eqnarray}
Recall that any positive solution of (\ref{3-1}) is an even function $u(-x)=u(x)$, satisfying $xu'(x)<0$ for $x \ne 0$ so that $\max _{(-1,1)} u(x)=u(0)$, and that any non-trivial solution of (\ref{3-5}) is of one sign, so that we may assume that $w(x)>0$ for $x \in (-1,1)$, see e.g., P. Korman \cite{K2}, \cite{K}. 
\medskip

If $f'(u)>(p-1)\frac{f(u)}{u}>0$ for $u>0$, it is well known that any positive solution of (\ref{3-1}) is {\em non-singular}, i.e.,    the problem (\ref{3-5}) admits only the trivial solution $w(x) \equiv 0$. Now suppose that $f'(u)>(p-1)\frac{f(u)}{u}>0$ holds only for $u>\gamma$, for some $\gamma >0$. It turns out that positive solutions of (\ref{3-1}), with maximum value greater than $\gamma$ are still non-singular, provided that $\int _u^{\gamma} f(t) \, dt < 0$ for all $u \in (0,\gamma)$.
The main result is  stated next.
It is customary to denote $F(u)=\int_0^u f(t) \, dt$.
\begin{thm}\lbl{thm:10}
Assume that $f(u) \in C^1(\bar R_+)$, and for some $\gamma >0$ it satisfies
\beq
\lbl{4-1}
 f(\gamma)=0, \s \mbox{and $\; f(u)>0$ on $(\gamma,\infty)$, }
\eeq 
\beq
\lbl{4-a}
f'(u)>(p-1)\frac{f(u)}{u}, \s \s \mbox{for  $u>\gamma$} \,,
\eeq
\beq
\lbl{4-2}
F(\gamma)-F(u)=\int _u^{\gamma} f(t) \, dt < 0, \s \mbox{for $u \in (0,\gamma)$} \,.
\eeq
Then   any positive solution of (\ref{3-1}), satisfying
\beq
\lbl{4-3}
u(0)>\gamma, \s \mbox{and} \s u'(1)<0 \,,
\eeq
is  non-singular, which means that    the linearized problem (\ref{3-5}) admits only the trivial solution.
\end{thm}

In case $p=2$ this result was proved in P. Korman \cite{K1}, while for general $p>1$ a weaker result, requiring that $f(u)<0$ on $(0,\gamma)$, was given  in J. Cheng \cite{J}, see also P. Korman \cite{K2}, \cite{K} for a different proof, and a more detailed description of the solution curve. Other multiplicity results on $p$-Laplace equations include \cite{Af}, \cite{A}, \cite{W} and \cite{R}.
\medskip

\pf
Assume, on the contrary, that the problem (\ref{3-5}) admits  a non-trivial solution $w(x)>0$. Let $x_0 \in (0,1)$ denote the point where $u(x_0)=\gamma$. Define
\[
q(x)=(p-1)(1-x) \varphi(u'(x))+\varphi'(u'(x)) u(x) \,.
\]
We claim that 
\beq
\label{3-6}
q(x_0)<0 \,.
\eeq
Rewrite (using that $(p-1) \varphi(t)=t\varphi'(t)$)
\[
q(x)=\varphi'(u'(x)) \left[(1-x)u'(x)+u(x) \right] \,.
\] 
Since $\varphi'(t)>0$ for all $t \ne 0$, it suffices to show that the function $z(x) \equiv (1-x)u'(x)+u(x)<0$ satisfies $z(x_0)<0$.
Indeed,
\[
z(x_0)=\int _{x_0}^1  \left[u'(x_0)-u'(x) \right] \,dx <0\,,
\]
which implies   the desired inequality (\ref{3-6}), provided we can prove that 
\beq
\lbl{3-7}
u'(x_0)-u'(x)<0 \,, \s \mbox{for $x \in (x_0,1)$}\,.
\eeq
The ``energy" function $E(x)=\frac{p-1}{p} {|u'(x)|}^p+F(u(x))$ is  seen by differentiation to be a constant, so that $E(x)=E(x_0)$, or
\[
\frac{p-1}{p} {|u'(x)|}^p+F(u(x))=\frac{p-1}{p} {|u'(x_0)|}^p+F(\gamma)  \,, \s \mbox{for all $x $}\,.
\]
By the assumption (\ref{4-2}), it follows that 
\[
\frac{p-1}{p} \left[ {|u'(x)|}^p-{|u'(x_0)|}^p \right]=F(\gamma)-F(u(x))<0 \,, \s \mbox{for $x \in (x_0,1)$}\,,
\]
 justifying (\ref{3-7}), and then giving (\ref{3-6}).
\medskip

Next, we claim that 
\begin{equation}
\label{3-8}
(p-1)w(x_0)\varphi(u'(x_0))-u(x_0)w'(x_0) \varphi'(u'(x_0))>0 \,,
\end{equation}
which implies, in particular, that
\begin{equation}
\label{3-9}
w'(x_0)<0 \,.
\end{equation}
Indeed, by a direct computation, using (\ref{3-1}) and (\ref{3-5}),
\[
\left[(p-1)w(x)\varphi(u'(x))-u(x)w'(x) \varphi'(u'(x)) \right]'= \left[ f'(u) -(p-1) \frac{f(u)}{u} \right]uw\,.
\]
The quantity on the right is positive on $(0,x_0)$, in view of our condition (\ref{4-a}). Integration over $(0,x_0)$, gives (\ref{3-8}).
\medskip

We have for all $x \in [-1,1]$
\beq
\lbl{3-8a}
\varphi'(u') \left(u'w'-u''w \right) =constant=\varphi'(u'(1)) u'(1)w'(1)>0  \,,
\eeq
as follows by differentiation, and using the assumption $u'(1)<0$. Hence
\beq
\lbl{3-8aa}
u'(x)w'(x)-u''(x)w(x)>0\,, \s \mbox{for $x \in (x_0,1)$}\,.
\eeq

Since  $f(u(x_0))=0$, it follows from the equation (\ref{3-1}) that $u''(x_0)=0$.  Then (\ref{3-8a}) implies
\begin{eqnarray}
\label{3-10}
& \varphi'(u'(1)) u'(1)w'(1)=\varphi'(u'(x_0)) u'(x_0)w'(x_0) \\
& =(p-1) \varphi(u'(x_0)) w'(x_0) \,. \nonumber
\end{eqnarray}

We need the following function, motivated by M. Tang \cite{T} (which was introduced in P. Korman \cite{K2}, and used in Y. An et al  \cite{A})
\[
T(x)=x \left[ (p-1) \varphi(u'(x))w'(x)+ f(u(x))w(x) \right]-(p-1) \varphi(u'(x))w(x) \,.
\]
One verifies that
\begin{equation}
\label{3-11}
T'(x)=p  f(u(x))w(x) \,.
\end{equation}
Integrating (\ref{3-11}) over $(x_0,1)$, and using (\ref{4-2}) and  (\ref{3-8aa}), obtain 
\[
T(1)-T(x_0)=
p  \int_{x_0}^1 f(u(x))w(x) \, dx
\]
\[
=p\int _{x_0}^1  \left[F(u(x))-F(\gamma)\right]'\frac{w(x)}{u'(x)} \,dx
\]
\[
=-p\int _{x_0}^1  \left[F(u(x))-F(\gamma)\right]\frac{w'(x)u'(x)-w(x)u''(x)}{{u'}^2(x)} \, dx<0 \,,
\]
which implies that
\[
L \equiv (p-1)\varphi(u'(1))w'(1)- (p-1) x_0 \varphi(u'(x_0))w'(x_0)+(p-1) \varphi(u'(x_0))w(x_0)<0 \,.
\]
On the other hand, using (\ref{3-10}), then (\ref{3-8}), followed by (\ref{3-9}) and (\ref{3-6}), we estimate the same quantity as follows
\begin{eqnarray} \nonumber
&  L>(p-1)\varphi(u'(x_0))w'(x_0)-(p-1)x_0 \varphi(u'(x_0))w'(x_0)+u(x_0)w'(x_0) \varphi'(u'(x_0)) \\\nonumber
& =w'(x_0)q(x_0)>0 \,,\nonumber
\end{eqnarray}
a contradiction.
\hfill$\diamondsuit$ \medskip

We remark that in case $f(0)<0$ it is possible to have a singular positive solution with $u'(1)=0$, so that the assumption $u'(1)<0$ is necessary.
\medskip

We now consider the problem
\begin{equation}
\label{3-1n}
 \varphi(u'(x))'+ \la f(u(x))=0 \s\s \mbox{for $-1<x<1$}, \s u(-1)=u(1)=0 \,,
\end{equation}
depending on a positive parameter $\la$.
The following result follows the same way as the Theorem $3.1$ in \cite{K2}.

\begin{figure}
\scalebox{0.9}{\includegraphics{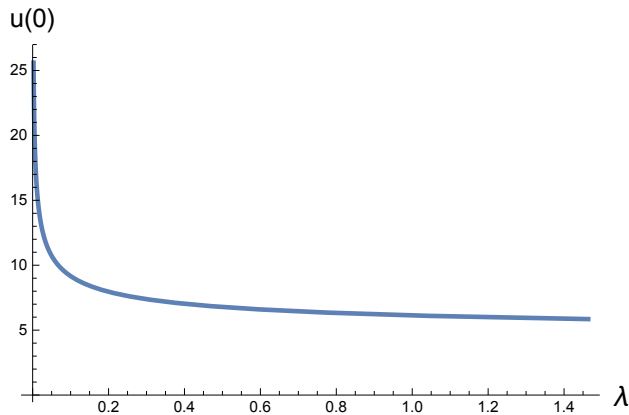}}
\caption{ The  curve of positive solutions for the problem (\ref{3-1n}), in case $p=3$ and $f(u)=u (u - 1) (u - 2) (u - 4)$.}
\label{fig}
\end{figure}

\begin{thm}\lbl{thm:2}
Assume that $f(u) \in C^1(\bar R_+)$, and the conditions (\ref{4-1}), (\ref{4-a}) and (\ref{4-2}) hold. Then there exists $0<\la _0 \leq \infty$ so that the problem (\ref{3-1n}) has a unique positive solution for $0<\la<\la _0$. All positive solutions, satisfying $u(0)>\gamma$,  lie on a continuous solution curve that is decreasing in the $(\la,u(0))$ plane (see Figure \ref{fig}). In case $f(0)<0$, one has $\la _0<\infty$, and at $\la =\la _0$ a positive solution with $u'(\pm 1)=0$ exists, and no positive solutions exist for $\la >\la _0$.
\end{thm}
\medskip

\noindent
{\bf Example} $\s$ In Figure \ref{fig} we present the solution curve of the problem (\ref{3-1n}) in case $p=3$ and $f(u)=u (u - 1) (u - 2) (u - 4)$. Here $\gamma=4$, and one verifies that the Theorem \ref{thm:2} applies. The {\em Mathematica} program to perform   numerical computations for this problem is explained in detail in \cite{KS} (it uses the shoot-and-scale method). The solution curve in Figure \ref{fig} exhausts the set of all positive solutions (since $\int _0^2 f(u) \, du<0$, there are no solutions with $u(0) =\max _{(-1,1)} u(x) \in (1,2)$).

\end{document}